\documentclass[12pt]{amsart}
\usepackage{amscd,amssymb}
\theoremstyle{latex 2e}
\newtheorem{thm}[subsection]{Theorem}
\newtheorem{lem}[subsection]{Lemma}
\newtheorem{prop}[subsection]{Proposition}

\newtheorem{con}[subsection]{Conjecture}

\newtheorem{remark}[subsection]{Remark}

\numberwithin{equation}{section}

\begin{document}

\title[Positive Quaternionic K\"ahler manifolds and symmetry rank: II] %
{Positive Quaternionic K\"ahler manifolds and symmetry rank: II}

\author{Fuquan Fang}
\address{Nankai Institute of Mathematics,
Weijin Road 94, Tianjin 300071, P.R.China}
\address{Department of Mathematics, Capital Normal University,
Beijing, P.R.China} \email{ffang@nankai.edu.cn}
\email{fuquanfang@eyou.com}

\thanks{**Supported by 973 project of Foundation Science of China, RFDP}

\begin{abstract}
Let $M$ be a positive quaternionic K\"ahler manifold of dimension $4m$.
If the isometry group $\text{Isom}(M)$ has rank at least
$\frac {m}2 +3$, then $M$ is isometric to $\Bbb HP^m$ or
$Gr_2(\Bbb C^{m+2})$. The lower bound for the rank is optimal if $m$ is even.
\end{abstract}
\maketitle

\section{Introduction}
\label{sec:intro}

\vskip2mm

\vskip4mm

A quaternionic K\"ahler manifold $M$ is an oriented Riemannian $4n$-manifold,
$n\ge 2$, whose holonomy group is contained in $Sp(n)Sp(1)\subset SO(4n)$.
If $n=1$ we add the condition that $M$ is Einstein and self dual.
Equivalently, there exists a $3$-dimensional subbundle $S$, of the
endmorphism bundle $\text{End}(TM, TM)$ locally
generated by three anti-commuting almost complex structures $I, J, K=IJ$
so that the Levi-Civita connection preserves $S$. It is well-known \cite{Be} that a
quaternionic K\"ahler manifold $M$ is always Einstein, and is necessarily
locally hyperK\"ahler if its Ricci tensor vanishes. A quaternionic K\"ahler
manifold $M$ is called {\it positive } if it has positive scalar curvature.
By \cite{Hi} (for $n=1$) and \cite{Sa} (for $n\ge 2$, compare \cite{Le} \cite{Le-Sa}) a positive
quaternionic K\"ahler manifold $M$ has a twistor space a complex Fano
manifold. Hitchin \cite{Hi} proved a positive quaternionic K\"ahler $4$-manifold
$M$ must be isometric to
$\Bbb CP^2$ or $S^4$. Hitchin's work was extended by Poon-Salamon \cite{PS} to
dimension $8$, which proves that a positive quaternionic K\"ahler $8$-manifold
$M$ must be isometric to $\Bbb HP^2$, $Gr_2(\Bbb C^{4})$ or $G_2/SO(4)$.

This leads to the Salamon-Lebrun conjecture:

{\it Every positive quaternionic
K\"ahler manifold is a quaternionic symmetric space.}

Very recently, the conjecture was further verified for $n=3$ in \cite{HH}, using
the approach initiated in \cite{Sa} \cite{PS} (compare \cite{Le-Sa}).
For a positive quaternionic K\"ahler manifold $M$, Salamon \cite{Sa} proved that
the dimension of its isometry group is equal to the index of certain twisted
Dirac operator, by the Atiyah-Singer index theorem, which is a characteristic
number of $M$ coupled with the Kraines $4$-form $\Omega $ (in analog with the
K\"ahler form), and it was applied to prove that the isometry group of $M$ is
large in lower dimensions (up to dimension $16$).

By \cite{Le-Sa} a positive quaternionic K\"ahler $4n$-manifold $M$ is simply
connected and the second homotopy group $\pi _2(M)$ is a finite group or
$\Bbb Z$, and $M$ is isometric to $\Bbb HP^n$ or $Gr_2(\Bbb C^{n+2})$
according to $\pi _2(M)=0$ or $\Bbb Z$.

An interesting question is to study positive quaternionic K\"ahler manifold
in terms of its isometry group. This approach
dates back to the work \cite{PS} for $n=2$ \cite{HH} for $n=3$ to proving the action is
transitive, and \cite{DS} \cite{PV} for cohomogeneity one actions (and hence the
isometry group must be very large). \cite{Bi} classified positive quaternionic
K\"ahler $4n$-manifolds with isometry rank $n+1$, using an approach on
hyper-K\"ahler quantizations. \cite{Fa} establishes a connectedness theorem and
using this tool the author proved that, a positive quaternionic K\"ahler
$4n$-manifolds of symmetry rank $\ge n-2$ must be either isometric to
$\Bbb HP^n$ or $Gr_2(\Bbb C^{n+2})$, if $n\ge 10$.

In this paper we will combine Morse theory of the momentum map on quaternionic K\"ahler
manifold developped in \cite{Ba} and the connectedness theorem in \cite{Fa} to prove the following

\vskip 2mm

\begin{thm}

Let $M$ be a positive quaternionic K\"ahler manifold of dimension $4m$.
Then the isometry group $\text{Isom}(M)$ has rank (denoted by
$\text{rank}(M)$) at most $(m+1)$, and $M$ is isometric to $\Bbb HP^m$ or
$Gr_2(\Bbb C^{m+2})$ if $\text{rank}(M) \ge \frac {m}2 +3$.
\end{thm}

\vskip2mm

Since the fixed point set of an isometric circle action is either a
quaternionic K\"ahler submanifold or a K\"ahler manifold. In the latter case
the fixed point set has dimension at most $2m$ (the middle dimension of the
manifold). Moreover, if a fixed point component is contained in
$\mu ^{-1}(0)$ then it must be a quaternionic K\"ahler submanifold, and if it
is in the complement $M-\mu ^{-1}(0)$ then it is K\"ahler (see \cite{Ba}).

\vskip 2mm

\begin{thm}

Let $M$ be a positive quaternionic K\"ahler $4m$-manifold $M$ with
an isometric $S^1$-action. If $N$ is a fixed point component of codimension
$4$. If $m\ge 3$ then $M$ is isometric to $\Bbb HP^m$ or $Gr_2(\Bbb C^{m+2})$.
\end{thm}

\vskip 2mm

The idea of proving Theorem 1.2 is as follows: by the above we know that
the fixed point component $N\subset \mu ^{-1}(0)$ if $4m-4\ge 2m+1$ (cf. \cite{Ba}
Remark 3.2). Furthermore, by \cite{GL} one knows that the quantization
$\mu ^{-1}(0)/S^1$ has dimension at most $4m-4$ (cf. \cite{DS}). We will prove in section 3 that
$N=\mu ^{-1}(0)$. Then we may use the equivariant
Morse equality to show that $M=\Bbb HP^m$ if $b_2(M)=0$ (cf. Lemma 4.1).
By \cite{Le-Sa} this implies easily Theorem 1.2.

With the help of Theorem 1.2, the proof of Theorem 1.1 follows by
using induction on the dimension and Theorems 2.1 and 2.4.

Theorem 1.1 is optimal if $m$ is even since the rank of $\widetilde{\text{Gr}}_4
(\Bbb R^{m+4})$ is $\frac m2 +2$. We conjecture that when $m$ is odd,
the lower bound for the rank in Theorem 1.1 may be improved by $1$, that is

\begin{con}
Let $M$ be a positive quaternionic K\"ahler manifold of dimension $8m+4$. Then
$M$ is isometric to $\Bbb HP^{2m+1}$ or $Gr_2(\Bbb C^{2m+3})$ if
$\text{rank}(M) \ge  m +3$.
\end{con}

\vskip 3mm

{\bf Acknowledgement:} The author is very grateful to F.
Battagalia for many invaluable discussions, who pointed out a gap
in a previous version, which leads to the discovery of Theorem
3.1.

\vskip10mm

\section{ Preliminaries}
\label{sec:prel}

\vskip2mm

\vskip4mm

In this section we recall some results on quaternionic K\"ahler
manifolds needed in later sections.

\vskip 2mm

Let $(M,g)$ be a quaternionic K\"ahler manifold of dimension $4n$. Let $F\to M$ be the
principal $Sp(n)Sp(1)$-bundle over $M$. Locally, $F\to M$ can be lifted to a principal
$Sp(n)\times Sp(1)$-bundle, i.e., the fiberwise double cover of $F$. Let $E$, $H$ be the
locally defined bundles
associated to the standard complex representation of $Sp(n)$ and $Sp(1)$
respectively. The complexified cotangent bundle $T^*M_\Bbb C$ is isomorphic to
$E\otimes _{\Bbb C} H$. The adjoint representations of $Sp(n)$ and $Sp(1)$
give two bundles $S^2E$ and $S^2H$ over $M$, respectively. Given the inclusion of
the holonomy algebra $sp(n)\oplus sp(1)$ into $so(4n)$, the bundle
$S^2E\oplus S^2H$
can be regarded as a subbundle of the bundle of $2$-forms $\Lambda ^2
T^*M_\Bbb C$. The bundle $S^2H$ has fiber the Lie algebra $sp(1)$ and the
local bases $\{I, J, K\}$ cooresponding to $i, j$ and $k\in sp(1)$
and consisting of three almost complex structures such that $IJ=-JI
=K$.

The Kraines $4$-form, $\Omega$, associated to a quaternionic K\"ahler
manifold $M$, is a non-degenerate closed form which is defined by
$$\Omega  =\omega _1\wedge \omega _1+\omega _2\wedge \omega _2+
\omega _3\wedge \omega _3$$
where $\omega _1$, $\omega _2$ and $\omega _3$ are the locally defined
$2$-forms associated to the almost complex structures $I, J$ and $K$.
The form $\Omega$ is globally defined and non-degenerate, namely $\Omega ^n$
is a constant non-zero multiple of the volume form. It is well-known that
$\Omega$ is parallel if and only if $M$ has holonomy in $Sp(n)Sp(1)$, if
$n\ge 2$. Moreover, by \cite{Sw} $M$ has holonomy in $Sp(n)Sp(1)$ if and only if
$\Omega $ is closed, provided $n\ge 3$.

For a quanternionic K\"ahler manifold, the almost complex structures may not be globally
defined, e.g., the quaternionic projective spaces $\Bbb HP^n$ does not admit an almost
complex structure. If $I, J, K$ are integrable and covariantly constant with respect the
metric, the holonomy group reduces to $Sp(n)$ and so the quaternionic K\"ahler manifold is
hyperk\"ahler. Wolf \cite{Wo} classified quaternionic symmetric spaces of compact type, they are
listed as $\Bbb HP^n$,  the complex Grassmannian $Gr_2(\Bbb C^{n+2})$, and the oriented real
Grassmannian $\widetilde {\text{Gr}}_4(\Bbb R^{n+4})$, and  exactly one quaternionic
symmetric space for each compact simple Lie algebra, namely $G_2/SO(4)$,
$F_4/Sp(3)Sp(1)$, $E_6/SU(6)Sp(1)$, $E_7/Spin(12)Sp(1)$, $E_8/E_7Sp(1)$.

As we mentioned in the introduction, so far quaternionic symmetric spaces are
the only known examples of positive quaternionic K\"ahler manifold.

\vskip 4mm

\begin{thm}[\cite{Le-Sa}]

(i) (Fininteness) For any $n$, there are, modulo isometries and rescalings,
only finitely many positive quaternionic K\"ahler $4n$-manifolds.

(ii) (Strong rigidity) Let $(M, g)$ be a positive quaternionic K\"ahler
$4n$-manifold. Then $M$ is simply connected and
$$\pi _2(M)=\{\begin{array} {lcc}
0, \text{   } (M,g)=\Bbb HP^n \\
\Bbb Z, \text{  } (M,g)=\text{Gr}_2(\Bbb C^{n+2})\\
\text{finite with $2$-torsion}, \text{  } \text{otherwise}
\end{array}$$

\end{thm}

\vskip2mm

A submanifold $N$ in a quaternionic K\"ahler manifold is called a {\it
quaternionic submanifold} if the quaternionic structure (i.e. locally
defined $I, J, K$) preserves the tangent bundle of $N$.

\vskip 2mm

\begin{prop}[\cite{Gr}]
Any quaternionic submanifold in a quaternionic K\"ahler manifold is totally
geodesic and quaternionic K\"ahlerian.
\end{prop}

\vskip2mm

\vskip2mm

\begin{thm}[\cite{Fa}]

Let $M$ be a positive quaternionic K\"ahler manifold of dimension $4m$.
Assume $f=(f_1,f_2): N\to M\times M$, where $N=N_1\times N_2$ and
$f_i: N_i\to M$ are quaternionic immersions of compact quaternionic K\"ahler
manifolds of dimensions $4n_i$, $i=1,2$. Let $\Delta $ be the diagonal of
$M\times M$. Set $n=n_1+n_2$. Then:

\noindent (2.3.1) If $n\ge m$, then $f^{-1}(\Delta)$ is nonempty.

\noindent (2.3.2) If $n\ge m+1$, then $f^{-1}(\Delta)$ is connected.

\noindent (2.3.3) If $f$ is an embedding, then for $i\le n-m$
there is a natural isomorphism, $\pi_i(N_1,N_1\cap N_2)\to
\pi_i(M,N_2)$  and a surjection for $i=n-m+1$.
\end{thm}

\vskip 2mm
As a direct corollary of (2.3.3) we have

\begin{thm}[\cite{Fa}]

Let $M$ be a positive quaternionic K\"ahler manifold of dimension $4m$.
If $N\subset M$ is a quaternionic K\"ahler submanifold of dimension $4n$,
then the inclusion $N\to M$ is $(2n-m+1)$-connected.
\end{thm}

\vskip10mm

\section{ HyperK\"ahler quotient and Quaternionic K\"ahler quotient}
\label{sec:hyper}

\vskip2mm

\vskip4mm

{\bf a. Hyperk\"ahler quotient}

\vskip 3mm

Let $M$ be a hyperk\"ahler manifold having a metric $g$ and
covariantly constant complex structures $I, J, K$ which behave
algebraically like quaternions:
$$I^2=J^2=K^2=-1; IJ=-JI=K$$

Let $G$ be a compact Lie group of isometries acting on $M$ and
preserving the structures $I, J, K$. The group $G$ preserves the
three K\"ahler forms $\omega _1, \omega _2, \omega _3$
corresponding to the three complex structures, so we may define
three moment maps $\mu _1, \mu _2, \mu_3$. These can be written as
a single map $$\mu : M\to   \mathfrak{g}^*\otimes \Bbb R^3$$

Let $$\mu _+=\mu _2+i\mu _3: M\to \mathfrak{g} ^* \otimes \Bbb C$$
where $\mathfrak{g}^*$ is the dual space of the Lie algebra of
$G$.

By \cite{HK} $\mu _+$ is holomorphic, and so $N=\mu _+^{-1}(0)$ is a
complex submanifold of $M$, with respect to the complex structure
$I$, therefore $N$ has an induced K\"ahler metric. By definition,
$\mu ^{-1}(0)=N\cap \mu _1^{-1}(0)$. The hyperk\"ahler quotient is
the quotient space $\mu ^{-1}(0)/G$, denoted by $M//G$. In
particular, if $\mu ^{-1}(0)$ is a manifold and the induced
$G$-action is free, then the hyperk\"ahler quotient $M//G$ is also
a hyperK\"ahler manifold. More generally, Dancer-Swann \cite{DS} proved
that the hyperk\"ahler quotient $M//G$ may be decomposed into the
union of hyperk\"ahler manifolds, according to the isotropy
decomposition of the $G$-action on $M$. However, it is not clear
at all if the decomposition of $M//G$ is a stratified topological
space satisfying the Goresky-MacPherson axioms, like in the
sympletic quotient case \cite{SL}.

In this section we will address to the structure of this
decomposition. For the sake of simplicity we consider only the
case of $G=S^1$ and the action is semi-free, i.e., free outside
the fixed point set.

Let us start with the standard example of isometric $S^1$-action
on quaternionic linear space $\Bbb H^n$ defined by
$$\varphi _t(u)=e^{2\pi it}u;  \hspace{1cm} t\in [0, 1)$$
where $i$ is one of the quaternionic units. With global
quaternionic coordinates $\{u^\alpha \}$, $\alpha =1, \cdots ,n$,
the standard flat metric on $\Bbb H^n$ may be written as:
$$ ds^2=\sum _\alpha d\bar u ^\alpha \otimes du^\alpha$$
where $\bar u^\alpha$ is the quaternionic conjugate of $u^\alpha$.

The Killing vector field $X$ of the above action is $\Bbb
H$-valued: $$X^\alpha (u)=iu^\alpha$$ which is triholomorphic.

Consider the $\Bbb H$-valued $2$-form
$$\omega=\sum _\alpha d\bar u ^\alpha \wedge du^\alpha$$
Observe that $\omega $ is purely imaginary since $\omega +\bar
\omega =0$. Note that $\omega =\omega _1 i+ \omega _2j+ \omega
_3k$, where $\omega _i$ is as above.

It is easy to see that the moment map (cf. \cite{Ga}) $$\mu ^X= \sum
_\alpha \bar u^\alpha iu^\alpha $$

Write $\Bbb H ^n=\Bbb C^n \oplus j \Bbb C^n$. The zero set
$M_0:=\mu ^{-1}(0)$ may be identified with the real algebraic variety of dimension $(4n-3)$:
$$\{(a, b)\in \Bbb C^n \oplus \Bbb C^n: |a|^2=|b|^2, \langle a , \bar b\rangle
=\sum _\alpha a^\alpha b^\alpha =0\}$$ The quotient $M_0/S^1$ is now an open cone over a
$(4n-5)$-dimensional manifold $W$:
$$W=\{(a, b)\in \Bbb C^n \oplus \Bbb C^n: |a|^2=|b|^2=1, \langle a , \bar b\rangle
=\sum _\alpha a^\alpha b^\alpha =0\}/S^1$$
In particular, if $n=1$, $M_0=\{0\}$ is a single point.

More generally, one may verify that, for any semi-free $S^1\subset Sp(n)$-action on $\Bbb H^n$ with
$\{ 0\}$ the only fixed point, the hyperk\"ahler quotient $\Bbb H^n//S^1$ is a topological cone
over a $(4n-5)$-dimensional real algebraic variety.

The following theorem is an analog of the
Sjamaar-Lerman theorem \cite{SL} (compare \cite{Gu}) for circle action in the hyperk\"ahler
case.

\vskip 3mm

\begin{thm}
Let $M^{4m}$ be a hyperk\"ahler manifold with an isometric
semi-free $S^1$-action preserving the hyperk\"ahler structure. If
$Y_1\subset M^{S^1}$ is a connected component of codimension $4n$
of the fixed point set in $\mu ^{-1}(0)$. Then $Y_1\subset M//S^1$
has an open neighborhood $U$ in $M//S^1$ which is diffeomorphic to
a fiber bundle over $Y_1$ with typical fiber a cone over a $(4n-5)$-dimensional
algebraic variety so that
$Y_1$ corresponds to the zero section. In particular, $Y_1$ is a
connected component of $M//S^1$ if $n=1$.
\end{thm}

\begin{proof} It suffices to prove that, for any given point $x\in Y_1$, there is an $S^1$-invariant
open neighborhood $U$ of $x\in M^{4m}$ such that $(U\cap \mu ^{-1}(0))/S^1$ is a fiber bundle over
$Y_1\cap U/S^1$ with typical fiber the hyperk\"ahler quotient $\Bbb H^n//S^1$.

Recall that $\omega =\omega _2+i\omega _3$ defines a complex symplectic structure on $M$. The usual
proof of the equivariant Darboux theorem for real sympletcic manifolds applies equally well in the
complex case; thus, there exists an $S^1$-equivariant open ball centered at $x\in M^{4m}$, such that
the complex symplectic manifold $(U, \omega )$ is equivaraintly complex symplectomorphic to
$(\Bbb H^m, i\sum _{\alpha =1}^m dz^\alpha \wedge dz^{\alpha +m})$, where $x$ corresponds to
the zero in $\Bbb H^m$. Identify $U$ with $\Bbb H^m$. By assumption, the fixed point component
$Y_1$ has codimension $n$,
without loss of the generality, we may assume that the quaternionic linear $S^1$-action on
$\Bbb H^m$ has a fixed point set the subspace $0\times \Bbb H^{m-n}$, and the action on $\Bbb
H^n\times 0$ is given by a semi-free $S^1\subset Sp(n)$-action with $\{0\}$ the only fixed
point. The moment map $\mu _+: U\to \Bbb C$ may be identified with the moment map of the linear
$S^1\subset Sp(m)$-action. Therefore $\mu _+^{-1}(0)$ is a complex algebraic variety, in fact
$\mu _+^{-1}(0)=\bar \mu _+^{-1}(0)\times \Bbb H^{m-n}$, where $\bar \mu _+: \Bbb H^m\to \Bbb
C$ is the moment map of the semifree $S^1\subset Sp(n)$ action on $\Bbb H^n$, which is a
complex algebraic variety with an induced K\"ahler metric (e.g., in the example above, this is
exactly the complex variety $\{(a, b)\in \Bbb C^n \oplus \Bbb C^n: \langle a, \bar b\rangle =0
\}$.)

Consider the variety as a symplectic manifold (at the zero, it is
reducible), by the classical Darboux theorem each irreducible
component may be identified with a standard symplectic ball with a
semi-free linear $S^1$-action, with center the only fixed point.
Therefore, $\mu ^{-1}(0)=\mu _1^{-1}(0)\cap \mu _+^{-1}(0)$ may be
identified with $(\bar \mu _1^{-1}(0)\cap \bar \mu
_+^{-1}(0))\times \Bbb H^{m-n}$, where $\bar \mu_1$ is the moment
map of the $S^1$-action on $\bar\mu _+^{-1}(0)$, which may be
identified (locally) with the zero set of the  moment map of an
semi-free linear $S^1\subset U(2n-1)$ on $\Bbb C^{2n-1}$ at every
irreducible component, which is a cone over a $(4n-4)$-dimensional
variety (cf. \cite{Gu}) Therefore, the quotient $(U\cap \mu
^{-1}(0))/S^1$ fibers over $Y_1\cap U/S^1$ with fiber a cone over
a $(4n-5)$-dimensional variety. Clearly, if $n=1$, then the cone
reduces to a single point. This proves the desired result.
\end{proof}

\vskip 2mm
\begin{remark}
It seems that the neighborhood $U$ may be chosen so that it is a fiber bundle over $Y_1$ with
typical fiber the hyperk\"ahler quotient $H^n//S^1$, where $S^1$ acts on $\Bbb H^n$ by the
isotropy representation at the fixed point set $Y_1$.
\end{remark}

\vskip 2mm
\begin{remark}
The same argument in the above proof extends trivially to more general situation. We will
come back to this point in some future paper.
\end{remark}

\vskip 4mm

{\bf b. Quaternionic K\"ahler quotient}

\vskip 3mm

Let $M$ be a quaternionic K\"ahler manifold with non-zero scalar curvature. If $G$ acts on $M$
by isometries, there is a well-defined moment map, which is a section
$\mu\in \Gamma (S^2H\otimes \mathfrak{g}^*)$ solving the equation
$$\langle \nabla \mu , X\rangle =\sum _{i=1}^3 \overline{I_iX}\otimes I_i$$
for each $X\in \mathfrak{g}$; where $\bar X=g(X, \cdot )$ denote the
$1$-form dual to $X$ with respect to the Riemannian metric. Equivalently, the above equation
may be written in the following form similar to the symplectic case
$$d\mu (X)=i_X\Omega$$

A nontrivial feature for quaternionic quotient is, the section $\mu$ is uniquely determined
if the scalar curvature is nonzero. Moreover, only the preimage of the zero section of the
moment map, $\mu ^{-1}(0)$, is well-defined.

\vskip2mm

\begin{thm}[\cite{GL}]
Let $M^{4n}$ be a quaternionic K\"ahler manifold with nonzero
scalar curvature acted on isometrically by $S^1$. If $S^1$ acts
freely on $\mu ^{-1}(0)$ then $\mu ^{-1}(0)/S^1$ is a quaternionic
K\"ahler manifold of dimension $4(n-1)$.
\end{thm}

Since the proof of the Galicki-Lawson's theorem is local, so if the circle acts freely on a
piece of the manifold is free, the same result applies well to the moment map on this piece.

\vskip 2mm

Recall that a Morse function $f$ is called {\it equivariantly
perfect} over $\Bbb Q$ if the equivariant Morse equialities hold,
that is if
$$\hat P_t(M)=\hat P_t(\mu ^{-1 }(0))+\sum t^{\lambda _F}\hat P_t(F)$$
where the sum ranges over the set of connected components of the
fixed point set, $\lambda _F$ is the index of $F$, and $\hat P_t$
is the equivariant Poincar\'e polynomial for the equivariant
cohomology with coefficients in $\Bbb Q$.

\vskip 2mm

\begin{thm}[\cite{Ba}]
Let $M^{4n}$ be a quaternionic K\"ahler manifold acted on
isometrically by $S^1$. Then the nondegenerate Morse function
$\|\mu \|^2$ is equivariantly perfect over $\Bbb Q$.
\end{thm}

\vskip 2mm

Let $f=\|\mu \|^2$. By \cite{Ba} the critical set of $f$ is the union
of the zero set $f^{-1}(0)=\mu ^{-1}(0)$ and the fixed point set
of the circle action.  Moreover, the zero set $\mu ^{-1}(0)$ is
connected, and a fixed point component is either contained in $\mu
^{-1}(0)$ or does not interesect with $\mu ^{-1}(0)$. The folllowing result is important for
this paper.

\vskip 2mm

\begin{prop}[\cite{Ba}]
Let $M^{4n}$ be a positive quaternionic K\"ahler manifold acted on
isometrically by $S^1$. Then every connected component of the
fixed point set, not contained in $\mu ^{-1}(0)$, is a K\"ahler
submanifold of $M- \mu ^{-1}(0)$ of real dimension less than or
equal to $2n$ whose Morse index is at least $2n$, with respect to the function $f$.
\end{prop}

For each quaternionic K\"ahler manifold $M$ with non-zero scalar curvature, following
\cite{Sw}, let $\mathfrak{u}(M)$ denote the $H^*/\{\pm 1\}$-bundle over $M$:
$$
\mathfrak{u}(M)=F\times _{Sp(n)Sp(1)}(H^*/\{\pm 1\})
$$
where $F$ is the principal $Sp(n)Sp(1)$-bundle over $M$. Let $\pi
: \mathfrak{u}(M)\to M$ denote the bundle projection. Obviously,
if $G$ is acts on $M$ by isometries, $G$ can be lifted to a
$G$-action on $\mathfrak{u}(M)$.  It is proved in \cite{Sw} that,
if the scalar curvature is positive, $\mathfrak{u}(M)$ has a
hyperk\"ahler structure which is preserved by the lifted
$G$-action. Moreover, the moment map $\mu $ is just the projection
of the moment map $\hat \mu$ of the lifted action on
$\mathfrak{u}(M)$ by $\pi$ (cf. \cite{Sw} Lemma 4.4).

\vskip 2mm

\begin{lem} Let $M$ be a positive quaternionic K\"ahler manifold
of dimension $4n$. Assume that $S^1$ acts on $M$ effectively by
isometries. Let $\mu \in \Gamma (S^2H)$ be its moment map. If
$N\subset \mu ^{-1}(0)$ is a fixed point component of codimension
$4$, then $N=\mu ^{-1}(0)$.
\end{lem}
\begin{proof} Let $\mathfrak{u}(M)$ be as above. By Proposition 4.2 of
\cite{DS}, at the fixed point $x\in N$, the isotropy
representation of $S^1$ in $SO(3)\cong
\mbox{Aut}(\mathfrak{u}(M)_x)$ is a finite group, where
$\mbox{Aut}(\mathfrak{u}(M)_x)$ is the isomorphism group of the
fiber at $x$ preserving the quaternionic structure. Therefore, the
preimage $\pi ^{-1}(N)$ is also a fixed point component of the
lifted $S^1$-action on $\mathfrak{u}(M)$, which has codimension
$4$.

By \cite{Sw} Lemma 4.4 we see that $\pi ^{-1}(N)\subset \hat \mu
^{-1}(0)$, where $\hat \mu$ is the moment map for the lifted
$S^1$-action on $\mathfrak{u}(M)$. Now $S^1$ acts on the normal
slice of $\pi ^{-1}(N)$ in $\mathfrak{u}(M)$ through a
representation in $Sp(1)$. For dimension reasoning, this
representation is faithful, otherwise, a finite order subgroup of
$S^1$ acts trivially on the whole manifold $\mathfrak{u}(M)$ and
so on $M$, a contradiction to the effectiveness of the action from
our assumption. Therefore, $S^1$ acts semi-freely on a
neighborhood of $\pi ^{-1}(N)$ in $\mathfrak{u} (M)$. By now we
may apply Theorem 3.1 to show that $\pi ^{-1}(N)$ is a connected
component of $\hat \mu ^{-1}(0)$. Since the moment map $\hat \mu$
projects to the moment map $\mu$, therefore $N$ is also a
connected component of $\mu ^{-1}(0)$. By \cite{Ba} $\mu ^{-1}(0)$
is connected, thus $N=\mu ^{-1}(0)$, the desired result follows.
\end{proof}

\vskip10mm

\section{Proof of Theorem 1.2}
\label{sec:proof 1.2}

\vskip2mm

\vskip4mm

Theorem 1.2 follows readily from the following Lemma and Theorem
2.1, where the dimension bound $m\ge 3$ implies that the fixed
point component of codimension $4$ has to be contained in $\mu
^{-1}(0)$, by Proposition 3.6.

\begin{lem}
Let $M$ be a positive quaternionic K\"ahler $4n$-manifold with an isometric
$S^1$-action where $n\ge 3$. Let $\mu$ be the moment map. Assuming $b_2(M)=0$.
If $N\subset \mu ^{-1}(0)$ is a fixed point component of dimension $4n-4$ of
the circle action, then $M$ is isometric to $\Bbb HP^n$.
\end{lem}
\begin{proof}[Proof]
By Lemma 3.7 $\mu ^{-1}(0)=N$, therefore $S^1$ acts trivially on
$\mu ^{-1}(0)$.

By Theorem 3.5
$$
\hat P_t(M)=\hat P_t(N)+\sum _Ft^{\lambda _F}\hat P_t(F)
$$
where $F$ runs over fixed point components outside $N$, and
$\lambda _F$ the Morse index of $F$. By Proposition 3.6 the Morse
index $\lambda _F\ge 2n$ and are all even numbers. Thus the
inclusion $N\to M$ is a $(2n-1)$-equivalence.

By \cite{Ba} Lemma 2.2  $\hat P_t(M)=P_t(M)P_t(BS^1)$. Since $S^1$
acts trivially on $F$ and $N$, we get that $\hat
P_t(F)=P_t(F)P_t(BS^1)$ and $\hat P_t(N)=P_t(N)P_t(BS^1)$. The
above identity reduces to
$$P_t(M)-P_t(N)=  \sum _Ft^{\lambda _F}P_t(F)$$
Observe that the last two terms of the left hand side is
$b_2(M)t^{4n-2}+t^{4n}$.

If $F$ is a fixed point component outside $\mu ^{-1}(0)$ such that
$\mbox{dim}_\Bbb RF>0$, we claim that $\mbox{dim}_\Bbb R F+\lambda
_F\le 4n-4$. Otherwise, by the above identity $\mbox{dim}_\Bbb R
F+\lambda _F=4n-2$ is impossible, and if the even integer
$\mbox{dim}_\Bbb R F+\lambda _F=4n$, we conclude that the
coefficients of $t^{4n-2}$ of the right hand side is also
non-zero, since $F$ must be a compact K\"ahler manifold (by
Proposition 3.6) and so $P_t(F)$ has nonzero coefficient at every
even degree not larger than the dimension.

Clearly the identity also shows that no isolated fixed point
outside $\mu ^{-1}(0)$ with Morse index $4n-2$, and there must
exist an isolated fixed point with Morse index $4n$.

Put these together, by Morse theory it follows that, up to
homotopy equivalence, $$M\simeq
 N\cup _{i} e^{\lambda _i}\cup e^{4n}$$ where $2n\le \lambda
 _i\le \mbox{dim}_\Bbb R F+\lambda _F\le 4n-4$, and $e^i$ denotes cell of dimension $i$.
Therefore $H^{4n-2}(M,N)=0$. By duality  $H_{2}(M-N)\cong
H^{4n-2}(M,N)=0$. Since the codimension of $N$ is $4$, it follows
that $H_2(M)\cong H_2(M-N)=0$. Therefore by Theorem 2.1 $M=\Bbb
HP^n$. The desired result follows.
\end{proof}

\vskip10mm

\section{ Proof of Theorem 1.1}
\label{sec:proof 2}

\vskip2mm

\vskip4mm

Let $M$ be a positive quaternionic K\"ahler manifold of dimension
$4m$. We call the rank of the isometry group $\text{Isom}(M)$ the
symmetry rank of $M$, denoted by $\text{rank}(M)$. By \cite{Fa} we
know that $\text{rank}(M)\le m+1$.

\vskip 2mm

\begin{proof}[Proof of Theorem 1.1]

Let $r=\text{rank}(M)$. Consider the isometric $T^r$-action on
$M$. Note that the $T^r$-action on $M$ must have non-empty fixed
point set since the Euler characteristic $\chi (M)>0$ by
\cite{Sa}. Consider the isotropy representation of $T^r$ at a
fixed point $x\in M$, which must be a representation through the
local linear holonomy $Sp(m)Sp(1)$ representation at $T_xM\cong
\Bbb H^m$. If there is a stratum (a fixed point set of an isotropy
group of rank $\ge 1$) of codimension $4$, then it must be
contained in $\mu ^{-1}(0)$ if $m\ge 3$ (by \cite{Ba} or
Proposition 3.6). By Theorem 2.1 and Lemma 4.1 the desired result
follows. Thus we can assume that at $x$, the isotropy
representation does not have any codimension $4$ linear subspace
fixed by some rank $1$ subgroup of $T^r$. Let $N$ be a maximal
dimensional submanifold of $M$ passing through $x$ fixed by a
circle subgroup of $T^r$.

Case (i): If $m=0(\text{mod }2)$.

By the above assumption $4m-8\ge \text{dim}N \ge 2m+4$ since
$\text{rank}(N)=r-1\ge \frac m2 +2$, by Lemma 2.1 of \cite{Fa}.
Note that $N$ is a quaternionic K\"ahler manifold since $N\subset
\mu ^{-1}(0)$. By Theorem 2.4 we see that $\pi _2(N)\cong \pi
_2(M)$. By Theorem 2.1 it suffices to prove $\pi _2(N)=0$ or $\Bbb
Z$. By induction we may consider $T^r$-action on $N$, and applying
Lemma 4.1 once again. Finally it suffices to consider the case
where a $16$-dimensional quaternionic K\"ahler submanifold of $M$,
$M^{16}$, with an effective isometric action by torus of
$\text{rank}\ge 5$. In this case there is a quaternionic K\"ahler
submanifold $M^{12}\subset M^{16}$ fixed by a circle group (cf.
\cite{Fa}). By Lemma 4.1, Theorem 2.1 and Theorem 2.4 $M^{16}=\Bbb
HP^4$ or $\text{Gr}_2 (\Bbb C^6)$, the desired result follows.

Case (ii): If $m=1(\text{mod }2)$.

Similar to the above $\text{dim}N \ge 2m+6$ for the same
reasoning. By Theorem 2.4 $\pi _2(N)\cong \pi _2(M)$. The same
argument by induction reduces the problem to the case of a
quaternionic K\"ahler submanifold of dimension $20$, $M^{20}$,
with an effective isometric torus action of $\text{rank}\ge 6$.
Once again the argument in \cite{Fa} shows that $M^{20}$ has a
quatertnionic submanifold $M^{16}$ of rank $\ge 5$. By (i) we see
that $M^{16}=\Bbb HP^4$ or $\text{Gr}_2 (\Bbb C^6)$. By Theorem
2.1 and Theorem 2.4 again we complete the proof.
\end{proof}


\end{document}